\documentstyle [twoside,12pt,amsmath,amssymb,amsthm]{article}
\newtheorem{thm}{Theorem}[section]
\newtheorem{cor}[thm]{Corollary}
\newtheorem{lem}[thm]{Lemma}
\newtheorem{prop}[thm]{Proposition}

\newcommand{\thmref}[1]{Theorem~\ref{#1}}
\newcommand{\lemref}[1]{Lemma~\ref{#1}}
\newcommand{\propref}[1]{Proposition~\ref{#1}}
\newcommand{\corref}[1]{Corollary~\ref{#1}}
\newcommand{\secref}[1]{Section~\ref{#1}}

\newcommand{\bx}{\hfill$\Box$\vspace{.6cm}}

\numberwithin{equation}{section}

\renewcommand\a{\alpha}         

\newcommand\g{\gamma}
\renewcommand\d{\delta}
\newcommand\e{\epsilon}
\renewcommand\l{\lambda}
\renewcommand\L{\Lambda}

\newcommand\G{\Gamma}
\newcommand\f{\frac}
\newcommand{\Z}{{\mathbb{Z}}}

\newcommand{\C}{{\mathbb{C}}}

\newcommand{\Q}{{\mathbb{Q}}}

\renewcommand\Re{\mbox{Re~}}

\renewcommand\O{{\cal O}}

\renewcommand\({\left(}         
\renewcommand\){\right)}

\font\cyrillic=wncyr10
\newcommand\Sha{\hbox{\cyrillic X}}
\newcommand{\can}{\hbox{can}}
\newcommand{\dd}{B}
\def\ST{\hbox{Shafarevich-Tate}}
\def\lr{\longrightarrow}
\def\can{\hbox{can}}

\begin{document}

\begin{center}
\vspace{-2cm} {\large \sc Non-vanishing of the Central Derivative
of Canonical Hecke L-functions}

 Stephen D. Miller and Tonghai Yang.

February 29, 2000

\end{center}
\sf

\setcounter{section}{1}

Let $K={\Q}(\sqrt{-D})$ be an imaginary quadratic field of
discriminant $-D<-4$, ${\O}$ its ring of integers, and $h$ its
ideal class number. A Hecke character
 $\chi$ of $K$ of conductor ${\frak f}$ is a called
``canonical" (\cite{Roh1}) if
\begin{eqnarray}
 & & \chi(\bar{\frak a})=\overline{\chi({\frak a})}
\mbox{ for each ideal } {\frak a} \mbox{ relatively prime to }
{\frak f}. \label{char1}
\\
& & \chi(\alpha{\O}) = \pm \alpha\mbox{ for principal ideals }
\alpha{\O} \mbox{ relatively prime to }{\frak f}. \label{char2}
\\
& & \mbox{The conductor }{\frak f} \mbox{ is divisible only by
primes dividing }D.  \label{char3}
\end{eqnarray}
Every Hecke character of $K$ satisfying (\ref{char1}) and
(\ref{char2}) is actually a quadratic twist of a canonical Hecke
character (see Section 2 for a precise description of these
characters and which fields have them).

 Let $L(s,\chi)$ denote the Hecke L-function of $\chi$, and $\L(s,\chi)$ its
completion; $\L(s,\chi)$ satisfies the functional equation
$\L(s,\chi)=W(\chi)\L(2-s,\chi)$, where $W(\chi)=\pm 1$ is the
root number. If $\chi$ is a canonical Hecke character with
$W(\chi)=1$, then
the central value $\L(1,\chi)\neq 0$
by a theorem of
 Montgomery and Rohrlich \cite{MR}.
 Of course, it automatically
vanishes when $W(\chi)=-1$ by the functional equation. The main
result of this paper is

\begin{thm}\label{mainthm}
Let $\chi$ be a canonical Hecke character whose root number
$W(\chi)=-1$.  Then the central derivative $\L'(1,\chi)\neq 0$.
\end{thm}
\vspace{.1cm}
 In \thmref{twists} we also prove that $\L'(1,\chi)\neq 0 $ when
$\chi$ is a small quadratic twist of a canonical character with
$W(\chi)=-1$.

When $D=p$ is a prime, canonical Hecke characters are closely
connected with the elliptic curves $A(p)$ extensively studied by
Gross \cite{Gro}. These curves are defined over
$F={\Q}(j(\f{1+\sqrt{-p}}{2} ) )$, where $j$ is the usual modular
$j$-function, and have complex multiplication by $\cal O$.
Combining \thmref{mainthm} and the above result of \cite{MR} with
Gross-Zagier \cite{GZ} and Kolyvagin-Logachev \cite{KL}, one has

\begin{cor}\label{apcor} Let $p>3$ be a prime congruent to 3 modulo 4.
Then

(a) \quad The Mordell-Weil rank of $A(p)$ is $$\hbox{rank}_{\Bbb
Z}A(p)(F)=\left\{\begin{array}{ll} h, & p\equiv 3~(\hbox{mod }8)
\\0, &p\equiv 7~(\hbox{mod}~8).
\end{array}     \right.$$

(b) \quad The Shafarevich-Tate group $\Sha (A(p)/F) $ is finite.
\end{cor}
In \cite{Gro}, Gross proved part (a) when $p\equiv 7~(\hbox{mod }
8)$ using a 2-descent.

In the next section we will outline the proof of \thmref{mainthm} and an
analog for quadratic twists (\thmref{twists}).  Sections 3, 4, and 5 are
devoted to analytic estimates used in the proofs of the theorems.  We
conclude in Section 6 with the proof of \corref{apcor} and other
arithmetic applications.

\subsection*{Acknowledgements}
 For very
large $D$, \thmref{mainthm} was originally obtained by D. Rohrlich
\cite{Roh4}. We are indebted to him for sharing his method with
us, as well as for allowing us to mention several of his results
here. Also, we thank him for his inspiration and for suggesting
this problem to us.  We are grateful to M. Baker, B. Conrad, N.
Elkies, B. Gross, K. Rubin, Z. Rudnick, D. Zagier, and S. Zhang
for discussions,
 and to Harvard University for their hospitality.  S.M. was
partially supported by an NSF post-doctoral fellowship as well as a Yale
Hellman
fellowship.  T.Y. was partially supported by an AMS Centennial fellowship
and NSF
grant DMS-9700777.

\section{Notation and Strategy}\label{notstrat}

We first recall some facts about canonical Hecke characters from
\cite{Roh2}.  They exist if and only if $D\equiv3~(\hbox{mod}~4)$
or is a multiple of 8.  Multiplying a canonical character by an
ideal class character always yields another canonical character.
This operation preserves the root number and defines natural
families of canonical Hecke characters.  When
$D\equiv3~(\hbox{mod}~4)$, there is exactly one family and it has
root number $\(\f2D\)$; when $D$ is a multiple of 8, there are two
families -- one has root number 1 and the other has root number -1.

To avoid confusion, we will sometimes write $\chi_{\can}$ for a
canonical Hecke character of $K$.  In this paper we consider Hecke
characters $\chi$ of $K$ satisfying conditions (\ref{char1}) and
(\ref{char2}), which are always of the form
\begin{equation}
\chi_{D,d}=\chi_{\can}\cdot(\e_d \circ N_{K/\Q}).
\label{chiDd}
\end{equation}
Here $d$ is a fundamental discriminant and
$\epsilon_d=(\frac{d}{})$ is the quadratic dirichlet character
with conductor $d$, prime to $D$.  The root number $W(\chi)$ is
explicitly computed in \cite{Roh2}.  In particular, when $D$ is
odd
\begin{equation}
W(\chi_{D,d})=\(\f2D\)\hbox{sign}(d).
\label{rootnum}
\end{equation}
From now on we will assume that $W(\chi)=-1$.
Set
\begin{equation}
B=\sqrt{DN{\frak f}}=\left\{\begin{array}{ll} D|d|, & D\mbox{ odd}
\\2D|d|, & 8\mid D.
\end{array}             \right.
\label{Bdef}
\end{equation}
The Hecke L-function is defined as
\begin{eqnarray}
& L(s,\chi)& =\sum_{{\frak a} \hbox{ integral}} \chi({\frak a}) (N
{\frak a})^{-s}\nonumber \\ & & =\sum_{\hbox{ideal classes }C}
L(s,\chi,C), \nonumber
\end{eqnarray}
where $L(s,\chi,C)$ is the partial L-series summed over integral
ideals in $C$.  Their completed L-functions are defined by
$$\L(s,\chi)=\(\f{B}{2\pi}\)^s \G(s)L(s,\chi)$$ and
$$\L(s,\chi,C)=\(\f{B}{2\pi}\)^s \G(s)L(s,\chi,C).$$

\begin{lem}\label{partl}
When $W(\chi)=-1$, $\L'(1,\chi)=0$ if and only if
$\L'(1,\chi,C)=0$ for each ideal class $C$ of $K$.
\end{lem}
{\bf Proof:} Associated to $\chi$ is a cuspidal new form $f$ of
weight 2 and level $B^2$ such that $L(s,f)=L(s,\chi)$.  So
Corollary 2 of \cite{GZ} implies that $\L'(1,\chi)=0$ if and only
if $\L'(1,\chi^\sigma)=0$ for every $\sigma \in
\hbox{Gal}(\bar{\Q}/{\Q})$.  On the other hand, by Theorem 1 of
\cite{Roh3}, $$ \{\chi^{\sigma}: \sigma \in \hbox{Gal}(\bar{\Q}/K)
\}
 =\{\chi \phi: \phi \hbox{ is an ideal class character of } K
 \},$$
 and $L(s,\chi)=L(s,\bar{\chi})$ by (\ref{char1}).  Thus
 $\L'(1,\chi)=0$ if and only if $\L'(1,\chi\phi)=0$ for all ideal
 class characters $\phi$ of $K$.
The ideal class characters are linearly
independent and $$
L(s, \chi\phi)= \sum_C \phi(C)L(s, \chi, C),
$$ so the lemma follows.\bx

To prove $\L'(1,\chi)\neq 0$, it now suffices to show
$\L'(1,\chi,c_1)\neq 0$ for the trivial class $c_1$ (i.e. the
class of principal ideals).  Since the root number $W(\chi)=-1$,
we have the functional equation
\begin{equation}
\L(s,\chi,c_1)=-\L(2-s,\chi,c_1).\label{funequ}\end{equation} By
Cauchy's theorem $$\L'(1,\chi,c_1)=\f{1}{2\pi
i}\(\int_{2-i\infty}^{2+i\infty}\L(s,\chi,c_1)\f{ds}{(s-1)^2}
-\int_{-i\infty}^{+i\infty}\L(s,\chi,c_1)\f{ds}{(s-1)^2} \).$$
Applying (\ref{funequ}) we arrive at the formula
\begin{equation}
\f12 \L'(1,\chi,c_1) = \f{1}{2\pi
i}\int_{2-i\infty}^{2+i\infty}\L(s,\chi,c_1)\f{ds}{(s-1)^2}.
\label{contour}
\end{equation}

It is clear from property (\ref{char2}) of $\chi$ that there is a
quadratic character  $\epsilon$ of $(\O/\frak f)^*$ such that
\begin{equation}
\chi(\alpha \O) = \epsilon(\alpha) \alpha. \label{epsdef}
\end{equation}   We
can express $L(s,\chi,c_1)$ as a sum over real and complex ideals:
\begin{equation}
L(s,\chi,c_1)=\sum_{n=1}^\infty \epsilon(n) n^{1-2s}
+\sum_{n=1}^\infty a_n n^{-s}, \label{randc}
\end{equation}
where
\begin{equation}
a_n=\sum_{\stackrel{N{{\frak a}} =n,~\frak a \ne \bar{{\frak
a}}}{\stackrel{\hbox{principal,}}{\hbox{integral}}}}
     \chi({\frak a})=
\sum_{\stackrel{u^2+Dv^2=4n}{u,v>0}}\e\(\f{u+\sqrt{-D}v}{2}\)u.
\label{andef}
\end{equation}

 Let \begin{equation}
  f(x)=\f{\G(0,x)}{x}=\f{1}{x}
\int_x^\infty e^{-t}\f{dt}{t}
 \label{fdef}\end{equation}
be the inverse Mellin transform of $\f{\G(s)}{(s-1)^2}$.  Indeed
\begin{eqnarray}
\int_0^\infty f(x)x^s\f{dx}{x}&=\int_0^\infty \int_x^\infty
x^{s-1}e^{-t}\f{dt}{t}\f{dx}{x}& \nonumber \\
&=\int_0^\infty\int_0^t \f{dx}{x} x^{s-1} e^{-t} \f{dt}{t}&
\nonumber \\ &=\f{1}{s-1}\int_0^\infty
t^{s-1}e^{-t}\f{dt}{t}=&\f{\G(s-1)}{(s-1)}=\f{\G(s)}{(s-1)^2},
\nonumber
\end{eqnarray}
so \begin{equation}f(y)=\f{1}{2\pi i}\int_{Re(s)=2}
\f{\G(s)}{(s-1)^2}y^{-s} ds .\label{mellinv}\end{equation}
Combining (\ref{contour}), (\ref{randc}), and (\ref{mellinv}) we
obtain
\begin{equation}
\f{1}{2}\L'(1,\chi,c_1) = \overbrace{\sum_{n=1}^\infty \epsilon(n)
n\cdot f(2\pi n^2/B )}^R+\overbrace{\sum_{n=1}^\infty a_n
f\(\f{2\pi n}{\dd} \)}^C. \label{rplusc}
\end{equation}
Formula (\ref{rplusc}) is essentially due to Rohrlich
({\cite{Roh4}), except that he expressed $R$ in terms of dirichlet
L-functions.

\subsection*{Examples: $D=8$ and $11$}

We will now illustrate (\ref{rplusc}) with the first two
discriminants which occur. Since both ${\Q}(\sqrt{-8})$ and
${\Q}(\sqrt{-11})$ have class number 1,
$$\L'(1,\chi,c_1)=\L'(1,\chi).$$ In order to compute it using
(\ref{rplusc}), we must first describe the character $\e:\({\cal
O}/{\frak f}\)^*\rightarrow \{\pm 1\}$. When $D=8$, ${\frak
f}=2\sqrt{-D}{\cal O}={\Z}8\oplus{\Z}\sqrt{-32}$, and $ \({\cal
O}/{\frak f}\)^*$ is generated by $\({\Z}/8\)^*$ and
$1+\sqrt{-2}$.  The character $\e(n)$ must restrict to
$\(\f{-8}{n}\)$ for $n\in \({\Z}/8\)^*$, and is thus determined by
its value on $1+\sqrt{-2}$.  In fact, $W(\chi)=\e(1+\sqrt{-2})$,
so in our case the values of $\e$ on the relatively-prime residue
classes are given in the following chart: \vspace{1cm}
\begin{center}
\begin{tabular}{c|c|rrrr}
 $\e(u+v\sqrt{-2})$ & $u$ &  1 &  3 &  5 &  7 \\
 \hline
         $v$          &     &    &    &    &    \\
\hline
         0          &     &  1 &  1 & -1 & -1 \\
         1          &     & -1 & -1 &  1 &  1 \\
         2          &     & -1 & -1 &  1 &  1 \\
         3          &     & -1 & -1 &  1 &  1 \\
\end{tabular}
\end{center}
\vspace{.5 cm}

When $D=11$, $\e\(\f{u+\sqrt{-11}v}{2}\)=\(\f{2u}{11}\)$.  We now
compute $L'(1,\chi)$ for these canonical characters, and check
(\ref{rplusc}) by comparing the known values of $L'(1,E)$ for the
associated elliptic curves $E$.

\begin{tabular}{l|cc}
        &   $D=8$ & $D=11$ \\
        \hline\hline
 Term $R$ with $n^2 \le 50$                                   &             1.82582357875147              &                               0.81497705252487                               \\
\hline
 Term $C$ with $n\le 50$                                      &             -0.28596530872740             &                             -0.0600975766040368                              \\
\hline
 $L'(1,\chi)=\(\f{2\pi}{B}\)\G(1)\L'(1,\chi)$ &             1.209401857169272             &                              0.862372296690396                               \\
        $\hspace{1.4cm}=\f{4\pi}{B}(R+C)$    &       &   \\
\hline

& & \\ Associated curve $E$ & $y^2=x^3+4x^2+2x$  &
 ~~~~~~~~$y^2+y=x^3-x^2-7x+10$
 \\   &
(\cite{Cr}, curve 256A)   & (\cite{Cr}, curve 121B)
\\\hline
 $L'(1,E)$ from \cite{Cr}
                  &               1.2094018572                &                                 .8623722967                                  \\
\hline\hline
\end{tabular}


\subsection*{Proof of \thmref{mainthm}}
By \lemref{partl} and (\ref{rplusc}), it suffices to prove $R
>|C|$.
 In the next section, we will prove that $R$ is bounded below by
\begin{eqnarray}
 R &>&\sum_{n=1}^\infty \l(n) n\cdot f\(\f{2\pi
n^2}{B}\) \nonumber \\
& > &  .5235B-.8458 B^{3/4}- .3951 B^{1/2}.\label{rlowerbd}
\end{eqnarray}
 Here $\l(n)$ is
Liouville's function -- the completely multiplicative function
which is -1 at each prime.

In \secref{trivbdsec}  we consider the special case $d=1$, and
bound term $C$ by \propref{Ctrivbd}:
 \begin{equation}
|C|<\left\{\begin{array}{ll} .2369 D, & D\mbox{ even}
\\.0269 D, &D\mbox{ odd}.
\end{array}             \right.
\label{thetrivbds}
\end{equation}

Having collected these estimates, proving $R > |C|$ is a simple
calculation. Indeed, if $D\ge 24$ is even, then $B=2D$, and $$R>
.5235 (2D) - .8458 (2D) (48)^{-1/4} - .3951 (2D) (48)^{-1/2}
=.2902 D,$$ so $$R>.2369 D>|C|.$$

If $D\ge 19$ is odd then
$$R>.5235 D - .8458\cdot D \cdot 19^{-1/4} - .3951\cdot D \cdot 19^{-1/2}
>.0277 D,$$
$$R>.0269 D >|C|.$$

There are only two values of $D$ not covered by this argument:
$D=8$ and $11$, which were dealt with in the examples.

\bx
%

\subsection*{Quadratic twists}
To prove non-vanishing for quadratic twists of canonical
characters, the bound (\ref{thetrivbds}) is not useful. In
\secref{rohrbdsec} we apply Rohrlich's method to obtain the
following bound on $C$ for $\chi_{D,d}$
(\propref{Crohrbd}):\footnote{The notation $A\ll B$ means
$A=O(B)$, i.e. there exists a positive constant $C$ such that
$|A|\le CB$.}
\begin{equation}|C|\ll D^{15/16+\d} |d|^{51/16+\d},\label{cddbd}
\end{equation}
where $\d>0$ is arbitrary and the implied constant depends only on
it.  Combining (\ref{rplusc}), (\ref{rlowerbd}), and (\ref{cddbd})
we conclude \vspace{.3cm}
\begin{thm}
\label{twists}
For any fixed $\d>0$, $$\L'(1,\chi_{D,d})\neq 0$$ for $|d| \ll
D^{1/35-\d}$ and $W(\chi_{D,d})=-1$.
\end{thm}
\vspace{.6cm}  {\bf Remarks.} When the root number $W(\chi)=1$,
similar non-vanishing results for twists were obtained in
\cite{Roh1}, \cite{RVY}, and \cite{Ya} for the central L-value.

 For canonical Hecke characters, Rohrlich
(\cite{Roh4}) computed $R$ as a contour integral of dirichlet
L-functions.  By shifting contours, $R$ can be expressed as the
sum of a residue and a remainder integral.  He showed the residue
is of size $\gg D$, and used Burgess' sub-convexity estimate to
bound the remainder integral by $\ll D^\a$, $\a<1$.  Also, he used
the method in Section 5 to show $C \ll D^\a$.  The power of the
main term is larger than that of the other two terms, and
positivity follows for large $D$.  However, the implied constants
one gets for these estimates are quite unfavorable.  In our proof
of \thmref{mainthm} we sacrifice the gain in the powers of $D$ for
a tie -- in favor of better constants.

\section{The Main Term $R$}

The purpose of this section is to prove (\ref{rlowerbd}).
We will show term $R$ is large and positive by eventually bounding
it from below by the following sum.

\begin{prop}\label{lampos}
$$\sum_{n=1}^\infty \l(n) n \cdot f\(\f{2 \pi n^2}{x}\)$$ is
always positive for $x>0$ and in fact
\begin{equation}
\sum_{n=1}^\infty \l(n) n \cdot f\(\f{2 \pi n^2}{x}\)> .5235 x -
.8458 x^{3/4}- .3951 x^{1/2} \label{lamlbd}
\end{equation} for $x>1$.

\end{prop}

{\bf Proof:} Using (\ref{mellinv}) and the identity
$$ \sum_{n=1}^\infty \l(n)n^{-s}=\f{\zeta(2s)}{\zeta(s)},$$ we
write
\begin{eqnarray}
&\sum_{n=1}^\infty \l(n) n \cdot f\(\f{2 \pi n^2}{x}\) & =
\f{1}{2\pi i}\int_{\Re s=2} \(\f{x}{2\pi}\)^s
\f{\G(s)}{(s-1)^2}\f{\zeta(4s-2)}{\zeta(2s-1)}ds \nonumber \\ & &
=\f{1}{2\pi i} \int_{\g}\(\f{x}{2\pi}\)^s
\f{\G(s)}{(s-1)^2}\f{\zeta(4s-2)}{\zeta(2s-1)}ds\nonumber \\ & &
+Res_{s=1} \(\f{x}{2\pi}\)^s
\f{\G(s)}{(s-1)^2}\f{\zeta(4s-2)}{\zeta(2s-1)} \nonumber \\ & & +
Res_{s=3/4}\(\f{x}{2\pi}\)^s
\f{\G(s)}{(s-1)^2}\f{\zeta(4s-2)}{\zeta(2s-1)}.\nonumber
\end{eqnarray}
 Here
$\g=C_1\cup C_2\cup C_3 \cup C_4\cup C_5$ is the contour
consisting of the union of the following five line segments:
\begin{eqnarray}
C_1 &\mbox{ from }1-i\infty&\mbox{ to }1-7i,\nonumber \\
 C_2&\mbox{ from
}1-7i&\mbox{ to }\f12-7i,\nonumber \\ C_3&\mbox{ from
}\f12-7i&\mbox{ to }\f12+7i,\nonumber \\ C_4&\mbox{ from
}\f12+7i&\mbox{ to }1+7i,\nonumber \\ C_5&\mbox{ from }1+7i&\mbox{
to }1+i\infty\nonumber\end{eqnarray} (7 is chosen because the
first critical zeroes of $\zeta(s)$ are approximately $\f12\pm
14.13472i$). The residue at $s=1$ is $\f{\pi x}{6}\approx
.523599x$ and the residue at $s=3/4$ is
$$Res_{s=3/4}\(\f{x}{2\pi}\)^s
\f{\G(s)}{(s-1)^2}\f{\zeta(4s-2)}{\zeta(2s-1)} = \f{2^{5/4}
x^{3/4} \G(3/4)}{\pi^{3/4}\zeta(1/2)}\approx -.845767 x^{3/4}.$$
 One can
easily estimate the integrals over $\g$ as follows.\footnote{All
computations were done using Mathematica v4.0 on an Intel Celeron
processor under Windows 98.}   First,
\begin{eqnarray}
&&\left| \int_{C_1}\(\f{x}{2\pi}\)^s
\f{\G(s)}{(s-1)^2}\f{\zeta(4s-2)}{\zeta(2s-1)} ds \right| \nonumber
\\
&&=\left|
\int_{C_5}\(\f{x}{2\pi}\)^s
\f{\G(s)}{(s-1)^2}\f{\zeta(4s-2)}{\zeta(2s-1)} ds      \right|
\nonumber
\\
& & \le \f{x}{2\pi} \int_{t=7}^\infty
\f{|\G(1+it)|}{t^2}\f{|\zeta(2+4it)|}{|\zeta(1+2it)|}dt \nonumber
\\ & & \le x (5\cdot10^{-7}).\nonumber \end{eqnarray} Next
\begin{eqnarray} &&\left|\int_{C_2} \(\f{x}{2\pi}\)^s
\f{\G(s)}{(s-1)^2}\f{\zeta(4s-2)}{\zeta(2s-1)} ds \right|\nonumber
\\ &&=\left|\int_{C_4} \(\f{x}{2\pi}\)^s
\f{\G(s)}{(s-1)^2}\f{\zeta(4s-2)}{\zeta(2s-1)} ds \right|\nonumber
\\ &&\le x \int_{\sigma=1/2}^1(2\pi)^{-\sigma}
\f{|\G(\sigma+7i)|}{(\sigma-1)^2+49}
\f{|\zeta(4\sigma-2+28i)|}{|\zeta(2\sigma-1+14i)|}d\sigma
\nonumber
\\ & &  \le x(2\cdot 10^{-6}).\nonumber \end{eqnarray}
 Finally,\begin{eqnarray}& &\left|\int_{C_3}\(\f{x}{2\pi}\)^s
\f{\G(s)}{(s-1)^2}\f{\zeta(4s-2)}{\zeta(2s-1)} ds \right|
\nonumber \\
& &
\le
\sqrt{\f{x}{2\pi}}\int_{t=-7}^{7}
\f{|\G(1/2+it)|}{1/4+t^2}\f{|\zeta(4it)|}{|\zeta(2it)|}dt
\nonumber \\
 &&\le
2.48218\sqrt{x}.\nonumber
 \end{eqnarray}
  Combining these estimates proves
(\ref{lamlbd}). For $x\ge 20,$
\begin{eqnarray}
& & .5235 x - .8458 x^{3/4} - .3951 x^{1/2} \nonumber  \\ & & \ge
(.5235- .8458\cdot 20^{-1/4} - .3951\cdot 20^{-1/2}  ) x \nonumber
\\ & & \ge .0351 x
> 0.\nonumber \end{eqnarray}
 The positivity for $x<20$ is
handled by the next lemma.\bx
\begin{lem}
\label{shortrangepos} For $0<x<20$, one has $$ f(\frac{2 \pi}x)>
\sum_{n=2}^\infty n \cdot f\(\f{2\pi n^2}{x}\). $$
\end{lem}
{\bf Proof:}
It is easy to see that for any $0< a<1$ and $ t > \frac{a}{1-a}$
$$ f(t) > a e^{-t}/t^2. $$ Take $a=\frac{\pi}{10+\pi} > .23$. Then
for $0< x < 20$, one has $$ f(\frac{2 \pi}x) > .23  \frac{x^2}{4
\pi^2} e^{-\frac{2 \pi}x}. $$ One the other hand, clearly, $f(x) <
e^{-x}/x^2$, and so
 \begin{eqnarray}
\sum_{n=2}^\infty n\cdot f\(\f{2\pi n^2}{x}\)
  \le & \sum_{n=2}^\infty n \f{x^2}{4\pi^2n^4}e^{-2\pi n^2/x}& \nonumber
\\
  \le &
\f{x^2}{4\pi^2}\sum_{n=2}^\infty n^{-3} e^{-2\pi n^2/x}.& \nonumber \\
\end{eqnarray}
 Since
$n^2\ge n+2$ for $n\ge 2$, this is
\begin{equation}
\le \f{x^2}{4 \pi^2} \f{e^{-8\pi/x}}{8} \sum_{n=0}^\infty e^{-2\pi
n/x}
=
\f{x^2
e^{-2\pi/x}}{4\pi^2}\frac{1}8\f{e^{-6\pi/x}}{(1-e^{-2\pi/x})}.
  \label{stepinshortrange}\end{equation}
Since $\frac{1}8\f{e^{-6\pi/x}}{(1-e^{-2\pi/x})}$ is clearly
increasing, it is thus bounded above in $0<x<20$ by its value
$\approx.181$ at $x=20$.  Therefore (\ref{stepinshortrange}) is
bounded above by $$ .19 \f{x^2}{4 \pi^2}e^{-2\pi/x} <
f(2\pi/x).$$\bx

\begin{prop}\label{othermults} (1)  If $m$ is any completely
multiplicative function with values $-1,0,$ or $1$, then

\begin{equation}\sum_{n=1}^\infty m(n) n \cdot f\(\f{2\pi n^2}{x}\)
 >0~,~x>0.\label{lowercharbd1}\end{equation}

(2)  If $ m_1$ and $ m_2$ are two distinct such functions with
$m_1(p) \ge m_2(p)$ for every prime p, then
\begin{equation}
\sum_{n=1}^\infty m_1(n) n\cdot f\(\f{2\pi n^2}{x}\)
>\sum_{n=1}^\infty m_2(n) n\cdot f\(\f{2\pi n^2}{x}\)
\label{lowercharbd2}
\end{equation}
for all $x>0$.
\end{prop}
{\bf Proof:}
We first assume
that the functions $m$, $m_1$, and $m_2$ differ from $\lambda$ at
only finitely many primes. Therefore, by \propref{lampos} and
induction, it suffices to prove (\ref{lowercharbd2}) under the following
conditions:

(a) \quad $m_2$ satisfies (\ref{lowercharbd1}) for all $x>0$.

(b) \quad $m_1$ and $m_2$ differ at exactly one prime, say  $p$,   and
    $m_1(p)=m_2(p) +1$.
Under assumption (b),  the difference is $$ \sum_{n=1}^\infty
m_1(n)n\cdot f\(\f{2\pi n^2}{x}\)-\sum_{n=1}^\infty
   m_2(n)n\cdot f\(\f{2\pi n^2}{x}\)$$
$$ =\sum_{p\mid n} (m_1(n)-m_2(n))n\cdot f\(\f{2\pi n^2}{x}\). $$
When $m_2(p)=-1$, $m_1(p)=0$, the difference is then $$
-\sum_{p\mid n} m_2(n)n\cdot f\(\f{2\pi n^2}{x}\)=
p\sum_{n=1}^{\infty} m_2(n)n\cdot f\(\f{2\pi n^2}{x/p^2}\) >0 $$
by assumption (a).  When $m_2(p)=0$, $m_1(p)=1$ and
$m_1(p^kn)=m_2(n)$ for $p\nmid n$. So the difference is $$
\sum_{p\mid n} m_1(n)n\cdot f\(\f{2\pi n^2}{x}\)
=\sum_{k=1}^{\infty}p^k
 \sum_{n=1}^{\infty} m_2(n)n\cdot f\(\f{2\pi n^2}{x/p^{2k}}\) >0
$$
by assumption (a) again.

In the general case, define $m^N(n)$ to be the completely
multiplicative function derived from $m$ by
\begin{equation}
m^N(p)
=
\left\{\begin{array}{ll} m(p), & p\le N
\\ \l(p) , & p>N.
\end{array}             \right.
\label{truncmul}
\end{equation}
Then $m^N$ differs from $\l$ at only a finite number of primes,
and thus $$\sum_{n=1}^\infty m^N(n)n\cdot f\(\f{2\pi n^2}{x}\)
\ge\sum_{n=1}^\infty \l(n)n\cdot f\(\f{2\pi n^2}{x}\) > 0.$$
Taking the limit as $N\rightarrow\infty$, $$\sum_{n=1}^\infty
m(n)n\cdot f\(\f{2\pi n^2}{x}\) \ge\sum_{n=1}^\infty \l(n)n\cdot
f\(\f{2\pi n^2}{x}\) > 0.$$  This completes the proof of part (1);
part (2) can be handled similarly.
 \bx

\begin{cor}\label{cor34} One has
$$ R \ge \sum_{n=1}^{\infty} \lambda(n) n \cdot f\(\frac{2 \pi
n^2}{\dd}\). \label{termIIbd} $$
\end{cor}

\vspace{.6cm} Combining \propref{lampos} with \corref{cor34}, one
obtains (\ref{rlowerbd}).

\section{The Trivial Bound on Remainder Term $C$}\label{trivbdsec}

In this section, we will only treat canonical characters, and prove
(\ref{thetrivbds}):

\begin{prop}
\label{Ctrivbd} When $d=1$ and $D\ge 7$, term $C$ is bounded by
\begin{equation} |C|< \left\{\begin{array}{ll} .0269 D, & D\mbox{
odd}
\\.2369 D, &D\mbox{ even}.
\end{array}             \right.
\label{Ctriveqn}\end{equation}
\end{prop}
{\bf Proof:} We first assume $D$ is odd, so $B=D$.  From
(\ref{andef}) we can bound $C$ term-wise, without appealing to
cancellation from the character.  To wit,

$$|C|\le \sum_{\stackrel{u,v>0}{u\equiv v~(\hbox{mod}~2)}} u
f\(\f{\pi}{2}(v^2+u^2/D)\).$$ Since $f(x)<e^{-x}/x^2$, $$|C| <
\sum_{\stackrel{u,v>0}{v\equiv u~(\hbox{mod}~2)}} u e^{-\f{\pi
u^2}{2 D}}\f{e^{-\pi v^2/2}}{\(\f{\pi}{2}(v^2+u^2/D)\)^2}$$ $$<
\sum_{u=1}^\infty u e^{-\f{\pi
u^2}{2D}}\f{4}{\pi^2}\sum_{\stackrel{v=1}{v\equiv
u~(\hbox{mod}~2)}}^\infty v^{-4}e^{-\pi v^2/2}.$$ The inside sum
is bounded by $$\sum_{v=1,~odd}^\infty v^{-4}e^{-\pi v^2/2}
\approx .20788,$$ so $$|C| \le .0843\sum_{u=1}^\infty u e^{-\f{\pi
u^2}{2D}}.$$

 If $D$ is
even, then $B=2D$,and the same argument shows that
\begin{eqnarray}|C|& <& 2\(\f{16}{\pi^2}\sum_{v=1}^\infty
v^{-4}e^{-\pi v^2/4} \)\(\sum_{u=1}^\infty u e^{-\pi u^2/D}\)
\nonumber \\& \le & 1.488 \sum_{u=1}^\infty u e^{-\pi u^2/D}.
\nonumber\end{eqnarray} The proof of (\ref{Ctriveqn}) now follows
from bound in the Lemma below.
 \bx

\begin{lem}
For $a\ge 1$,
\begin{equation}
\sum_{n=1}^\infty ne^{-n^2/a} < a/2. \label{poiss}
\end{equation}
\end{lem}
This conclusion actually holds for all $a>0$.

\vspace{.5cm}

 {\bf
Proof:} Using the Poisson summation formula applied to
$|n|e^{-n^2/a}$, $$\sum_{n=1}^\infty n e^{-n^2/a} = \int_0^\infty
ne^{-n^2/a}dn +2 \sum_{r=1}^\infty\int_0^\infty n
e^{-n^2/a}\cos(2\pi r n)dn.$$ The first integral is $a/2$ and the
others are actually negative.  This is because $$\int_0^\infty
ne^{-n^2/a}\cos(2\pi r n)dn = \f{a}{2}\(1-e^{-a\pi^2 r^2}2\pi
r\sqrt{a}\int_0^{\pi r\sqrt{a}}e^{t^2}dt\)$$ (cf. \cite{GR},
17.13.27) and $$\int_0^{\pi r\sqrt{a}}e^{t^2}dt > 1 +
\int_1^{a\pi^2 r^2}\f{e^t}{2\sqrt{t}}dt = 1+\f{1}{2\pi r
\sqrt{a}}\(e^{a\pi^2r^2}-e\)> \f{e^{a\pi^2r^2}}{2\pi r
\sqrt{a}}.$$\bx

\section{Rohrlich's Bound on Remainder Term $C$}\label{rohrbdsec}

We now return to the general case of a canonical character twisted
by $\e_d=\(\f{d}{}\)$.  The method here is adapted from
\cite{Roh1} and \cite{Roh4}.

\begin{prop}
\label{Crohrbd} For any $\d>0$, term $C$ is bounded by $$|C| \ll
D^{15/16+\d}|d|^{51/16+\d},$$ where the implied constant depends
only on $\d$.
\end{prop}

{\bf Proof:} Set $A(t)=\sum_{n<t}a_n$.  Integration by parts gives
\begin{equation}C=\int_{D/4}^\infty f\(\f{2\pi
t}{B}\)\f{d}{dt}A(t)dt = -\int_{D/4}^\infty
A(t)\f{d}{dt}f\(\f{2\pi
t}{B}\)dt,\label{rohrbdparts}\end{equation} because there are no
complex ideals of norm $<D/4$.  By \cite{Roh1}, p. 553(27), $A(t)$
is bounded above by $$|A(t)|\ll
t^{5/4}D^{-5/16+\d}|d|^{19/16+\d}$$ for $D>8$ and $t>0$ (the
implied constant again depends only on $\d$).  Along with the
inequalities $$0<-\f{d}{dt}f\(\f{2\pi t}{B}\)< \f{B}{2\pi t^2}
e^{-2\pi t/B}\(1+\f{B}{2\pi t}\),$$ (\ref{rohrbdparts}) implies
\begin{eqnarray}
|C|& \ll& D^{-5/16+\d}|d|^{19/16+\d}\int_{D/4}^\infty
\left[t^{5/4}\f{B}{2\pi t^2} e^{-2\pi t/B}\(\f43+\f{B}{2\pi
t}\)\right]dt\nonumber \\ &\ll&
D^{-5/16+\d}|d|^{19/16+\d}\int_{D/4}^\infty
\f{d}{dt}\left[-\f{B^2}{2\pi^2} t^{-3/4} e^{-2\pi t/B}
\right]dt\nonumber \\
  &  \ll& D^{-17/16+\d}|d|^{19/16+\d}B^2.\nonumber\end{eqnarray}
  Since either $B=D|d|$
    or $2D|d|$, this completes the proof.\bx

\section{Arithmetic Applications}\label{arithapp}

Having completed their proofs, we will now give some arithmetic
applications of Theorems~\ref{mainthm} and~\ref{twists}, including
\corref{apcor}.

Let $j$ be the $j$-invariant of a fixed isomorphism class of
elliptic curves with complex multiplication (CM) by $\cal O$. Then
$H=K(j)$ is the Hilbert class field of $K$.  We can extend any
Hecke character $\chi$ of $K$ to one on $H$ by $$ \psi= \chi \circ
N_{H/K}. $$ When $\chi$ satisfies (\ref{char1}) and (\ref{char2}),
$$ \psi({\frak A}^\sigma) = \psi({\frak
A})^\sigma,~~\sigma\in\hbox{Gal}(H/\Bbb Q) $$ for every ideal
$\frak A$ of $H$ relatively prime to the conductor of $\psi$. By
Theorem 9.1.3 and Lemma 11.1.1 of \cite{Gro}, there is a unique
elliptic $\Q$-curve $A$ over $H$ with $$j(A)=j\mbox{  and
}L(s,A/H)=L(s,\psi)L(s,\bar{\psi}).$$ (Here we recall that a
``$\Q$-curve'' is an elliptic curve over a number field which is
isogenous to all of its Galois conjugates.) Furthermore, $A$
descends to two isogenous elliptic curves over the subfield
$F={\Q}(j)$ (\cite{Gro}, Theorem 10.2.1).  By abuse of notation we
will also refer to these curves as $A$.  Let
$B=\hbox{Res}_{F/{\Q}} A$ be the abelian variety over $\Q$
obtained from $A$ by restriction of scalars.  When $D=p$ is prime,
Gross proved (\cite{Gro}, Theorem 15.2.5) that $T=\hbox{End}_K B
\otimes \Bbb Q$ is a CM number field of degree $2h$; thus $B$ is
also a CM abelian variety.  This result actually extends to
composite $D$ via a different argument:

\begin{lem}\label{lem61}
(a) Let $T$ be the subfield of $\C$ generated by $\chi(\frak a)$,
where $\frak a$ runs over all ideals of K prime to $\chi$'s
conductor.  Then $T$ is a CM number field of degree $2h$, and
$\Phi =\{ \sigma: T \rightarrow {\Bbb C}~|~ \sigma  \hbox{ trivial
on } K\}$ is a CM type of T.

(b) B is a CM abelian variety of type $(T, \Phi)$.
\end{lem}

{\bf Proof:} For each embedding $\sigma:T\rightarrow \C$ fixing $K$,
$\sigma\circ
\chi$ is
another canonical Hecke character of $K$, and thus it is of the form
$\chi\phi$, where $\phi$ is an ideal class character of $K$.  By Theorem 1
of \cite{Roh3}, $\sigma \mapsto \phi$  actually  gives a one-to-one
correspondence between the complex
embeddings of $T$ into $\Bbb C$ fixing $K$, and the ideal class characters
of $K$.
Thus $[T:K]=h$ and $[T:{\Q}]=2h$.  It is a general fact that $T$ is a CM
number field; in this case it can easily be verified using property
(\ref{char1}) of $\chi$.

By \cite{Sh}, Theorem 10, there is a CM abelian variety $B'/\Q$ of type
$(T,\Phi)$ associated to $\chi$, and it is unique up to isogeny.  In
particular
\begin{equation}
L(s, B') =\prod_{\sigma: T^+ \lr \Bbb C} L(s, \chi^{\sigma})
         =\prod_{\phi} L(s, \chi \phi),\label{equ61}
\end{equation} where $T^+$ is the maximal totally-real subfield of $T$.
On the other hand,
$$L(s,B)=L(s,A/F)=L(s,\psi)=\prod_{\phi}L(s,\chi\phi).$$ This
shows $L(s,B)=L(s,B')$, so a theorem of Faltings \cite{Fa}
guarantees $B$ and $B'$ are isogenous, proving (b).\bx

\begin{lem}
\label{lem62}
Let $\chi$ be a Hecke character of $K$ of the form (\ref{chiDd}).  Let $A$
be an associated $\Q$-curve over $F={\Q}(j)$ with $j$-invariant $j$, and
let $B=\hbox{Res}_{F/\Q}A$.  If $\hbox{ord}_{s=1}L(s,\chi)\le 1$ then

(a) The Mordell-Weil ranks of $A$ and $B$ are given by
$$\hbox{rank}_{{\Z}} A(F) = \hbox{rank}_{{\Z}} B({\Q}) = h\cdot
\hbox{ord}_{s=1}L(s,\chi).$$

(b) The Shafarevich-Tate groups $\Sha(A/F)$ and $\Sha(B/{\Q})$ are finite.
\end{lem}

{\bf Proof:} Since the Mordell-Weil and Shafarevich-Tate groups of
$A$ over $F$ are identical to those of $B$ over $\Q$, it is
sufficient to prove the Lemma for $B$.  Let $f$ be the normalized
weight 2 new-form associated to $\chi$ as in the proof of
\lemref{partl}. The field generated by $f$'s Fourier coefficients
is generated by $\chi(\frak a) + \overline{\chi(\frak a)}$, and is
thus $T^+$.  Equation (\ref{equ61}) implies
$$L(s,B)=\prod_{\sigma: T^+ \lr {\C}} L(s, f^{\sigma}). $$  Now
the Lemma follows from a result of Kolyvagin and Logachev
(\cite{KL}).\bx

Combining \lemref{lem62} with the non-vanishing theorems above
(Theorems \ref{mainthm} and~\ref{twists}) and in \cite{MR}, one
gets the following two corollaries.

\begin{cor}\label{cor63}
Let $\chi=\chi_{\hbox{can}}$ be a canonical Hecke character of $K$, and
let
$A$ and $B=\hbox{Res}_{F/{\Q}}A$ respectively be associated $\Q$-curves
and CM abelian
varieties as above.   Then

(a) The Mordell-Weil ranks of $A$ and $B$ are given in terms of the root
number $W(\chi)$ by
$$
\hbox{rank}_{\Z} A(F) = \hbox{rank}_{\Z} B({\Q}) =
\left\{\begin{array}{ll} h, & W(\chi)=-1
\\0, & W(\chi)=1.
\end{array}     \right.$$
In particular, when $D$ is odd, these ranks are $h$ or zero depending on
whether $D\equiv 3$ or 7$\mod 8$.

(b) The \ST~groups $\Sha(A/F)$ and $\Sha(B/{\Q})$ are finite.

 \end{cor}

{\bf Proof of \corref{apcor}:} Take $D=p$,
$j=j\(\f{1+\sqrt{-p}}{2}\)$, $A=A(p)$, and apply
\corref{cor63}.\bx

\begin{cor}\label{cor64}
Let $\chi=\chi_{D,d}$ be a Hecke character of $K$ of the form
(\ref{chiDd}). Let $A$ and $B=\hbox{Res}_{F/{\Q}} A$ be as above, and
fix any $\d>0$.
If $|d| \ll D^{1/35-\d}$ (the implied constant depending on $\d$) and
$W(\chi_{D,d})=-1$, then

(a) The Mordell-Weil ranks of $A$ and $B$ are
$$\hbox{rank}_{{\Z}} A(F) = \hbox{rank}_{{\Z}} B({\Q})=h.$$

(b) The \ST~groups $\Sha(A/F)$ and $\Sha(B/{{\Q}})$ are finite.
 \end{cor}

Finally, we wish to point out that when $D$ is prime, all $\Q$-curves over
$F$ are associated to Hecke characters of the form (\ref{chiDd}), though
this is not true for every composite $D$.  See \cite{Na} for
a more-precise description.

\vspace{1cm}

\begin{tabular}{lcl}
 Stephen D. Miller         &  & Tonghai Yang                 \\
 Department of Mathematics &  & Department of Mathematics    \\
 Yale University           &  & State University of New York \\
 P.O. Box 208283           &  & Stony Brook, NY 11794-3651   \\
 New Haven, CT 06520       &  & {\tt thyang@math.sunysb.edu}       \\
{\tt stephen.miller@yale.edu}   &  &
\\
\end{tabular}

\end{document}